\input amstex

\documentstyle{amsppt}
\document
\input xypic
\pageheight{45pc}
\NoBlackBoxes

\define\Q{\Bbb Q}

\define\Z{\Bbb Z}
\define\N{\Bbb N}

\define\C{\Bbb C}

\define\cal{\Cal}

\magnification 1200
\topmatter
\rightheadtext{On a Hasse principle}
\title On a Hasse principle for Mordell-Weil groups
\endtitle
\author Grzegorz Banaszak
\endauthor
\address
Department of Mathematics, Adam Mickiewicz University,
Pozna\'{n}, Poland
\endaddress
\email banaszak\@amu.edu.pl
\endemail

\abstract In this paper we establish a Hasse principle concerning 
the linear dependence over $\Z$ of nontorsion points in the 
Mordell-Weil group of an abelian variety over a number field.

\endabstract
\endtopmatter

\subhead  1. Introduction
\endsubhead

\medskip
Let $A$ be an abelian variety over a number 
field $F.$ Let $v$ be a prime of ${\cal O}_F$ 
and let $k_v := {\cal O}_F / v.$ Let $A_{v}$ denote the 
reduction of $A$ for a prime $v$ of good reduction and let 
$$r_v\, :\, A(F) \rightarrow A_v (k_v)$$
be the reduction map. Put ${\cal R} := End_F (A).$
Let $\Lambda$ be a subgroup of $A (F)$ and let 
$P \in A(F).$ A natural question arises whether the condition 
$r_{v}(P) \in r_v ( \Lambda)$ 
for almost all primes $v$ of ${\cal O}_F$ implies that 
$P\in \Lambda.$ This question was posed by W. Gajda in 2002.
The main result of this paper is the following theorem.

\proclaim{Theorem 1.1} Let $P_1, \dots, P_r$ be elements
of $A(F)$ linearly independent over ${\cal R}.$ 
Let $P$ be a point of $A(F)$ such that ${\cal R}\, P$ is a free
${\cal R}$ module. The following conditions are equivalent:
\roster
\item[1]\,\,$P\in \sum_{i = 1}^r \, \Z\, P_i$
\item[2]\,\,$r_{v}(P) \in \sum_{i = 1}^r \, \Z\, r_{v}(P_i)$ 
for almost all primes $v$ of ${\cal O}_F.$
\endroster
\endproclaim
\medskip

\noindent
In the case of the multiplicative group $F^{\times}$ 
the problem analogous to W. Gajda's question has already been solved 
by 1975. Namely, A. Schinzel, [Sch, Theorem 2, p. 398],
proved that for any $\gamma_1,\dots, \gamma_r \in F^{\times}$
and $\beta \in F^{\times}$ such that 
$\beta = \prod_{i = 1}^r \gamma_{i}^{n_{v,i}}\,\, \text{mod}\,\, v$
for some $n_{i,v}, \dots, n_{r, v} \in \Z$ for 
almost all primes $v$ of ${\cal O}_F$ 
there are $n_1,\dots, n_r \in \Z$ such that
$\beta = \prod_{i = 1}^r \gamma_{i}^{n_i}.$ Theorem of A. Schinzel
was proved again by Ch. Khare [Kh] 
using  methods of C. Corralez-Rodrig{\'a}{\~ n}ez and R. Schoof
[C-RS]. Ch. Khare used this theorem to prove that every family of one 
dimensional strictly compatible $l$-adic representations
comes from a Hecke character. 
\medskip

\noindent
Theorem 1.1 strengthens the results of [BGK2], [GG] and [We].
Namely T. Weston [We] obtained an analogue of
Theorem 1.1 with coefficients in $\Z$ 
for ${\cal R}$ commutative. T. Weston did not assume that 
$P_1,\dots, P_r$ is a basis over ${\cal R},$ however there was
some torsion ambiguity in the statement of his result.   
In [BGK2] together with W. Gajda and P. Kraso{\' n} 
we proved Theorem 1.1 for elliptic curves without CM and more 
generally for a class of
abelian varieties with $End_{{\overline F}} (A) = \Z.$
We also got a general result for all abelian varieties
[BGK2] Theorem 2.9 in the direction of Theorem 1.1. 
However in Theorem 2.9 loc. cit. the 
coefficients are in ${\cal R}$ and there is also a
coefficient in ${\N}$ associated with the point $P.$ 
Recently W. Gajda and K. G{\' o}rnisiewicz 
[GG] Theorem 5.1, strengthened [BGK2] Theorem 2.9 implementing some 
techniques of M. Larsen and R. Schoof [LS] and showing that
the coefficient associated with the point $P$ 
in [BGK2] Theorem 2.9 is equal to 1. The coefficients 
in [GG] Theorem 5.1 are still in ${\cal R}.$
Very recently A. Perucca has proven the Theorem 5.1 of [GG] 
(see [Pe] Corollary 5.2) using her $l$-adic support problem result, 
see loc. cit. At the end of this paper we prove 
Theorem 5.1 of [GG] by some of our methods from the proof of 
Theorem 1.1.
\medskip

Although not explicitely presented in our proofs, 
this paper makes an essential use of 
results on Kummer Theory for abelian varieties, originally developed 
by K. Ribet [Ri], and results of G. Faltings [Fa], J-P. Serre and J. 
Tate [ST], A. Weil [W], J. Zarhin [Za] and other important results about 
abelian varieties. The application of these results comes by quoting some 
results of [BGK1], [BGK2] and [Bar] where Kummer Theory and results of G. Faltings, 
J-P. Serre, J. Tate, A. Weil and J. Zarhin where key ingredients. 
This form of the exposition of our paper makes the proofs of our results 
concise and transparent for the reader.  
  
\subhead 2. Proof of Theorem 1.1
\endsubhead 
Let $L/F$ be an extension of number fields and let $w$ denote a 
prime ideal in ${\cal O}_L$ over a prime $v$ of good reduction, such that
$v$ does not divide $l.$ It follows by [BGK1], Lemma 2.13 
that the reduction map
$$r_w \, :\, A(L)_l \rightarrow A_w (k_w)_l$$
is injective, where $G_l$ denotes the $l$-torsion
part of an abelian group $G.$ 
The main ingredients in the proof of the 
Lemma 2.13 loc. cit. are [ST], Theorem 1 and 
Weil conjectures for abelian varieties, 
proven by Weil [W]. For additional information about the injectivity of 
the reduction map
$$r_w \, :\, A(L)_{tor} \rightarrow A_w (k_w)$$  
see also [K] p. 501-502. We will use several times 
in this paper the following result of S. Bara{\' n}czuk
which is a refinement of the Theorem 3.1 of [BGK2] and 
Proposition 2.2 of [BGK3].

\proclaim{Lemma 2.1} ([Bar], Th. 5.1)
Let $l$ be a prime number. 
Let $m_1, \dots, m_s \in \N \cup \{0\}$ and let 
$m := max \{m_1, \dots, m_s\}.$ Let $L/F$ be a finite extension
and let $Q_1, \dots, Q_s \in A(L)$ 
be independent over ${\cal R}.$ 
There is a family of primes $w$ of ${\cal O}_L$ 
of positive density such that
\, $r_w (Q_i)$ has order $l^{m_i}$ in $A_w (k_w)_l$ 
for all $1 \leq i \leq s,$
\endproclaim

\proclaim{Corollary 2.2}
Let $m \in \N.$ Let $Q_1, \dots, Q_s \in A(F)$ be independent over
${\cal R}$ and let $T_1, \dots T_s \in A[l^{m}].$  Let $L := F(A[l^m]).$ 
There is a family of primes $w$ of ${\cal O}_L$  
of positive density such that for the prime $v$ of ${\cal O}_F$
below $w:$  

\roster
\item[1]\,\, $r_w (T_1), \dots, r_w (T_s) \in  A_v (k_v),$
\item[2]\,\, $r_v (Q_i) = r_w (T_i)$ in $A_v (k_v)_l$ 
for all $1 \leq i \leq s.$
\endroster
\endproclaim

\demo{Proof} Observe that the points 
$Q_1 - T_1, \dots, Q_s - T_s$
are linearly independent over ${\cal R}$ in $A(L).$
Indeed if:
$$\sum_{i=1}^s \beta_i (Q_i - T_i) = 0 $$
for some $\beta_1, \dots, \beta_s  \in {\cal R}$, then 
$$\sum_{i=1}^s l^m \,\beta_i \, Q_i  = 0.$$ 
Hence $l^m\, \beta_i = 0$ for all
$1 \leq i \leq s.$ Since ${\cal R}$ is a torsion free 
abelian group it shows that $\beta_1 = \dots = \beta_s = 0.$ 
It follows by Lemma 2.1 that there is a family of primes $w$
of ${\cal O}_L$ of positive density such that 
$r_w (Q_i - T_i) = 0$ in $A_w (k_w)_l.$
Since $Q_1, \dots, Q_s \in A(F),$ it follows that
$r_w (Q_i - T_i) = r_w (Q_i) - r_w (T_i) = r_v (Q_i) - r_w (T_i)$
for the prime $v$ of ${\cal O}_F$ below $w.$ 
Hence  we get $r_w (T_i)  = r_v (Q_i) \in A_v (k_v)_l$ 
for all $1 \leq i \leq s.$ 
\enddemo

\demo{Proof Theorem 1.1} 
It is enough to prove that (2) implies (1). 
By Theorem 2.9 [BGK2] there is an $a \in \N$ and elements 
$\alpha_1, \dots, \alpha_r \in {\cal R}$ such that:
$$a P = \sum_{i=1}^r\, \alpha_i\, P_i. \tag{2.3}$$

\noindent
{\bf Step 1.} Assume that $\alpha_i \in \Z$ for all $1 \leq i \leq r.$
We will show (cf. the proof of Theorem 3.12 of [BGK2]) that 
$P\in \sum_{i = 1}^r \, \Z\, P_i.$
Let $l^k$ be the largest power of $l$
that divides $a.$ Lemma 2.1 shows that for any $1 \leq i \leq r$ there
are infinitely many primes $v$ such that 
$r_v(P_1) = \dots = r_v(P_{i-1}) = r_v(P_{i+1}) = \dots = r_v(P_r) = 0$ 
and $r_v(P_i)$ has order equal to $l^k$ in
$A_v(k_v)_l.$ By (2.3) we get $a r_v(P) = \alpha_i r_v(P_i)$ 
Moreover by assumption (2) of the theorem,
$r_v(P) = \beta_i r_v(P_i)$ 
for some $\beta_i \in \Z.$ Hence
$$(\alpha_i{-}a \beta_i) r_v(P_i) = 0$$ in
$A_v(k_v)_l.$ This implies that 
$l^k$ divides $\alpha_i$ for all $1 \leq i \leq r.$ So
By (2.3) we get 
$${a \over l^k} P = \sum_{i = 1}^r  {\alpha_{i} \over
l^k}P_{i} +  T ,\tag{2.4}$$
for some $T \in A(F)[l^k].$ Again, by Lemma 2.1 there are infinitely
many primes $v$ in ${\cal O}_F$ 
such that $r_v(P_i) =  0$ in $A_v (k_v)_l$ for all 
$1 \leq i \leq r.$ In addition 
$r_v (P) \in \sum_{i = 1}^r \, \Z\, r_v (P_i)$
for almost all $v.$ So (2.4) implies that $r_v(T) = 0,$ for
infinitely many primes $v.$ This contradicts the injectivity of $r_v,$ 
unless $T = 0.$ Hence
$${a \over l^k} P = \sum_{i = 1}^r {\alpha_{i} \over l^k} P_i. \tag{2.5}$$ 
Repeating  the above argument for primes dividing ${a
\over l^k}$ shows that condition (1) holds.
\medskip

\noindent
{\bf Step 2.}
Assume $\alpha_i \notin \Z$ for some $i.$ Observe that $\alpha_i$
is an endomorphism of the Riemann lattice ${\cal L},$  such that
$A(\C) \cong {\C}^g / {\cal L}.$ To make the notation simple, we will 
denote again by $\alpha_i$ the endomorphism $\alpha_i \otimes 1$ acting on 
$T_l (A) \cong {\cal L} \otimes \Z_l.$ 
Let $P(t) := \,\, \text{det} (t \text{Id}_{{\cal L}} - \alpha_{i})
\in \Z[t],$ be the characteristic polynomial of $\alpha_i$ acting on
${\cal L}.$ Let $K$ be the splitting 
field of $P(t)$ over $\Q.$ We take $l$ such that it splits in $K$ and 
$l$ does not divide primes of bad reduction. Since $P(t)$ has 
all roots on ${\cal O}_K$ and is also the 
characteristic polynomial of $\alpha_i$ on $T_l (A),$ we see that
$P(t)$ has all roots in $\Z_l$ by the assumption on $l.$
If $P(t)$ has at least two different roots in ${\cal O}_K,$ 
we easily find a vector $u \in T_l (A)$ which is not an eigenvector 
of $\alpha_i$ on $T_l (A).$ 
If $P(t)$ has a single root $\lambda \in {\cal O}_K$ then  
$P(t) = (t - \lambda )^{2g}$ and we must have 
$\lambda \in \Z$ because we are in characteristic 0.
Hence $P(t) = (t - \lambda )^{2g}$ is the characteristic polynomial of 
$\alpha_i$ as an endomorphism of ${\cal L}.$ Since $\alpha_i \notin \Z$ 
we find easily  $u \in {\cal L}$
such that $u$ is not an eigenvector of $\alpha_i$ acting on 
$T_l(A).$ In any case there is $u \in T_l(A)$ which is not an 
eigenvector of $\alpha_i$ acting on $T_l(A).$
Rescaling if necessary, we can 
assume that $u$ is not divisible by $l$ in $T_l (A).$ Hence for $m \in \N$
and $m$ big enough we can see that the coset $u + l^m T_l (A)$ 
is not an eigenvector of $\alpha_i$ acting on $T_l (A)/ l^m T_l (A).$ 
Indeed, if $\alpha_i u \equiv c_m u \,\, \text{mod}\,\, l^m T_l (A)$
for $c_m \in \Z/l^m$ for each $m \in \N,$ then $c_{m+1} u \equiv c_m u
\,\, \text{mod}\,\, l^m T_l (A).$ Because $u$ is not divisible by $l$ 
in $T_l (A),$ this implies that $c_{m+1} \equiv c_m \,\, \text{mod}\,\, l^m$ 
for each $m \in \N.$ 
But this contradicts the fact that $u$ is not an eigenvector
of $\alpha_i$ acting on $T_l (A).$
Consider the natural isomorphism of Galois and ${\cal R}$ 
modules $T_l (A)/ l^m T_l (A) \cong A[l^m].$ We put $T \in A[l^m]$
to be the image of the coset $u + l^m T_l (A)$ via this isomorphism.
Put $L := F(A[l^m]).$ By Corollary 2.2 we choose a prime $v$
below a prime $w$ of ${\cal O}_L$ such that  
\roster
\item"{(i)}" \,\, $r_w (T) \in A_v (k_v)_l,$
\item"{(ii)}" \,\, $r_v (P_j) = 0$ 
for all $j \not= i$ and $r_v (P_i) = r_w (T)$ in $A_v (k_v)_l.$ 
\endroster
From (2.3) and (ii) we get
$a r_v (P) = \alpha_i\, r_v (P_i) = \alpha_i\, r_w (T) \,\,\, 
\text{in}\,\,\,
A_v (k_v)_l.$ Hence for the prime $w$ in ${\cal O}_L$ over $v$ we get
in $A_w (k_w)_l$ the following equality: 
$$a r_w (P) = \alpha_i\, r_w (P_i) = \alpha_i\, r_w (T) \tag{2.4}$$
By assumption (2) and (ii) there is $d \in \Z,$ such that
$a r_v (P) = a d \, r_v (P_i) = a d\, r_w (T) 
\,\,\, \text{in}\,\,\, A_v (k_v)_l.$
Hence, for the prime $w$ in ${\cal O}_L$ over $v,$ we get
in $A_w (k_w)_l$ the following equality:  
$$a r_w (P) = a d \, r_w (P_i) = a d\, r_w (T). \tag{2.5}$$
Since $r_w$ is injective, the equalities (2.4) and (2.5) give: 
$$\alpha_i\, T = a d \, T \quad \text{in}\quad A[l^m].$$
But this contradicts the fact that $T$ is not an eigenvector
of $\alpha_i$ acting on $A[l^m].$ It proves that
$\alpha_i \in \Z$ for all $1 \leq i \leq r,$ but this case 
has been taken care of already in step 1 of this proof. $\qed$
\enddemo
\medskip

\noindent
\proclaim{Corollary 2.6} Let $A$ be a simple abelian variety. Let 
$P_1, \dots, P_r$ be elements
of $A(F)$ linearly independent over ${\cal R}.$ Let $P$ be a nontorsion 
point of $A(F).$ The following conditions are equivalent:
\roster
\item[1]\,\,$P\in \sum_{i = 1}^r \, \Z\, P_i$
\item[2]\,\,$r_{v}(P) \in \sum_{i = 1}^r \, \Z\, r_{v}(P_i)$ 
for almost all primes $v$ of ${\cal O}_F.$
\endroster
\endproclaim
\demo{Proof} This is an immediate consequence of Theorem 1.1.
Indeed, 
for a nontorsion point $P$ the ${\cal R}$-module ${\cal R}\, P$ is a free 
${\cal R}$-module since $D = {\cal R} \otimes_{\Z} \Q$ is a division algebra
because $A$ is simple. $\qed$
\enddemo
\medskip

\noindent
\proclaim{Corollary 2.7} Let $A$ be a simple abelian variety.
Let $P$ and  $Q$ be nontorsion elements of $A(F).$
The following conditions are equivalent:
\roster
\item[1]\,\,$P = m\, Q$ for some $m \in \Z,$
\item[2]\,\,$r_{v}(P) = m_v\, r_v (Q)$ 
for some $m_v \in \Z$ for almost all primes $v$ of ${\cal O}_F.$
\endroster
\endproclaim
\demo{Proof} This is an immediate consequence of Corollary 2.6 because 
${\cal R} P$ and ${\cal R} Q$ are free ${\cal R}$-modules since 
$A$ is simple. $\qed$
\enddemo
\medskip

\noindent
The following proposition is the Theorem 5.1 of [GG]. We give here a 
new proof of this theorem using some methods of the proof
of Theorem 1.1. 

\proclaim{Proposition 2.8} Let $A$ be an abelian variety 
over $F.$ Let $P_1, \dots, P_r$ be elements
of $A(F)$ linearly independent over ${\cal R}.$ 
Let $P$ be a point of $A(F)$ such that ${\cal R}\, P$ is a free
${\cal R}$ module. The following conditions are equivalent:
\roster
\item[1]\,\,$P\in \sum_{i = 1}^r \, {\cal R} \, P_i$
\item[2]\,\,$r_{v}(P) \in \sum_{i = 1}^r \, {\cal R} \, r_{v}(P_i)$ 
for almost all primes $v$ of ${\cal O}_F.$
\endroster
\endproclaim
\demo{Proof} 
Again we need to prove that (2) implies (1). Let us assume (2).
By [BGK2], Theorem 2.9 there is an $a \in \N$ and elements 
$\alpha_1, \dots, \alpha_r \in {\cal R}$ such that
equality (2.3) holds. Let $l$ be a prime number such that $l^k || a$
for some $k > 0.$ Put $L := F(A[l^k])$ and take arbitrary $T \in A[l^k].$
By Corollary 2.2 we can choose a prime $v$ below a prime $w$ of ${\cal O}_L$ 
such that  
\roster
\item"{(i)}" \,\, $r_w (T) \in A_v (k_v)_l,$
\item"{(ii)}" \,\, $r_v (P_j) = 0$ 
for all $j \not= i$ and $r_v (P_i) = r_w (T)$ in $A_v (k_v)_l.$ 
\endroster
From (2.3) and (ii) we get
$a r_v (P) = \alpha_i\, r_v (P_i) = \alpha_i\, r_w (T) \,\,\, 
\text{in}\,\,\,
A_v (k_v)_l.$ Hence we get the following equality in $A_w (k_w)_l:$ 
$$a r_w (P) = \alpha_i\, r_w (P_i) = \alpha_i\, r_w (T) \tag{2.9}$$
By assumption (2) and (ii) there is $\delta \in {\cal R},$ such that
$a r_v (P) = a \delta \, r_v (P_i) = a \delta \, r_w (T) = 0 
\,\,\, \text{in}\,\,\, A_v (k_v)_l.$
Hence we get the following equality in $A_w (k_w)_l$ :  
$$a r_w (P) = a \delta \, r_w (P_i) = a \delta \, r_w (T) = 0. \tag{2.10}$$
By injecivity of $r_w,$ the equalities (2.9) and (2.10) imply: 
$$\alpha_i\, T = 0 \quad \text{in}\quad A[l^k].$$
This shows that $\alpha_i$ maps to zero 
in $End_{G_F}\, (A[l^k]).$ By [Za] Corollary 5.4.5, cf. 
the proof of Lemma 2.2 of [BGK2], we have a natural isomorphism:
$${\cal R} / l^k {\cal R}\,\, \cong \,\, End_{G_F}\, (A[l^k]).$$
Hence we proved that $\alpha_i \in l^k {\cal R}$ for all $1 \leq i \leq r.$
So 
$${a \over l^k} P = \sum_{i = 1}^r  \beta_{i} P_{i} +  T^{\prime} ,\tag{2.11}$$
where $\beta_i \in {\cal R}$ for all $1 \leq i \leq r$ and 
$T^{\prime} \in A(F)[l^k].$ By Lemma 2.1 there are infinitely
many primes $v$ in ${\cal O}_F$ 
such that $r_v(P_i) =  0$ in $A_v (k_v)_l$ for all 
$1 \leq i \leq r.$ In addition 
$r_v (P) \in \sum_{i = 1}^r \, {\cal R} \, r_v (P_i)$
for almost all $v.$ So (2.11) implies that $r_v(T^{\prime}) = 0,$ for
infinitely many primes $v.$ Hence $T^{\prime} = 0$ by the injectivity
of $r_v$ (see [BGK1] Lemma 2.13). Hence
$${a \over l^k} P = \sum_{i = 1}^r \beta_{i} P_i. \tag{2.12}$$ 
Repeating  the above argument for primes dividing ${a
\over l^k}$ finishes the proof of the proposition.
$\qed$
\enddemo

\subhead 3. Remark on Mordell-Weil ${\cal R}$ systems 
\endsubhead 
Let ${\cal R}$ be a ring with identity. 
In the paper [BGK1] the Mordell-Weil 
${\cal R}$ systems have been defined. In [BGK2] we investigated 
Mordell-Weil ${\cal R}$ systems satisfying certain
natural axioms $A_1 - A_3$ and $B_1 - B_4.$ We also assumed
that ${\cal R}$ was a free
$\Z$-module. Let us consider 
Mordell-Weil ${\cal R}$ systems which are associated to
families of $l$-adic representations $\rho_l\, :\, G_F \rightarrow 
GL (T_l)$ such that $\rho_l (G_F)$ contains an open image of 
homotheties. Since Theorem 2.9 of [BGK2] and 
Theorem 5.1 of [Bar] were proven for Mordell-Weil ${\cal R}$ systems,
then Proposition 2.8 and its proof 
generalize for the Mordell-Weil ${\cal R}$ systems. This shows 
that Theorem 2.9 of [BGK2], which is stated for Mordell-Weil 
${\cal R}$ systems, holds with $a = 1.$ 
Let us also assume that there is a free $\Z$-module
${\cal L}$ such that ${\cal R} \subset End_{\Z} ({\cal L})$ and for 
each $l$ there
is an isomorphism ${\cal L} \otimes \Z_l \cong T_l$ such that the 
action of ${\cal R}$ on $T_l$ comes from its action on ${\cal L}.$ 
Abelian varieties are principal examples of Mordell-Weil
${\cal R}$ systems  satisfying all the requirements 
stated above with ${\cal R} = End_{F} (A).$  
Then Theorem 1.1 generalizes also for Mordell-Weil ${\cal R}$ systems
satisfying the above assumptions because we can apply again
Theorem 2.9 of [BGK2] and Theorem 5.1 of [Bar].  
\medskip

\noindent
{\it Acknowledgements}:\quad
The research was partially financed by 
the research grant N N201 1739 33 of the Polish Ministry of Science 
and Education and by the Marie Curie Research Training Network "Arithmetic 
Algebraic Geometry" MRTN-CT-2003-504917.

\enddocument